\documentclass[11pt]{article}

\usepackage{latexsym}
\usepackage{amssymb}

\newtheorem{thm}{Theorem}

\begin{document}
{
\begin{center}
{\Large\bf
The Nevanlinna-type formula for the matrix Hamburger moment problem in a general case.}
\end{center}
\begin{center}
{\bf S.M. Zagorodnyuk}
\end{center}

\section{Introduction.}
Recall that the matrix Hamburger moment problem consists of
finding a left-continuous non-decreasing matrix function $M(x) = ( m_{k,l}(x) )_{k,l=0}^{N-1}$
on $\mathbb{R}$, $M(-\infty)=0$, such that
\begin{equation}
\label{f1_1}
\int_\mathbb{R} x^n dM(x) = S_n,\qquad n\in \mathbb{Z}_+,
\end{equation}
where $\{ S_n \}_{n=0}^\infty$ is a prescribed sequence of Hermitian $(N\times N)$ complex matrices (moments),
$N\in \mathbb{N}$.
The  moment problem~(\ref{f1_1}) is said to be {\it determinate} if it has a unique solution and {\it indeterminate}
in the opposite case.

\noindent
This problem was introduced in~1949 by Krein~\cite{cit_100_K}, and he described all solutions
in the case when the corresponding J-matrix defines a symmetric operator with maximal defect numbers.
This result appeared without proof in~\cite{cit_200_K}
(Berezansky in 1965 proved the main fact in this theory of Krein: the convergence of the series
from the polynomials of the first kind, even for the operator moment problem~\cite[Ch.7, Section 2]{cit_700_Ber}).
Under similar conditions, some descriptions of solutions were obtained by
Kovalishina~\cite{cit_300_K}, by Lopez-Rodriguez~\cite{cit_400_L} and by Dyukarev~\cite{cit_500_D}.

In the scalar case, a description of all solutions of the moment problem~(\ref{f1_1}) can be found, e.g.,
in~\cite{cit_600_Akh},\cite{cit_700_Ber} for the nondegenerate case, and in~\cite{cit_800_AK} for
the degenerate case.

Set
\begin{equation}
\label{f1_3}
\Gamma_n = \left(
\begin{array}{cccc} S_0 & S_1 & \ldots & S_n\\
S_1 & S_2 & \ldots & S_{n+1}\\
\vdots & \vdots & \ddots & \vdots\\
S_n & S_{n+1} & \ldots & S_{2n}\end{array}
\right),\qquad n\in \mathbb{Z}_+.
\end{equation}
It is well known that the following condition
\begin{equation}
\label{f1_4}
\Gamma_n \geq 0,\qquad n\in \mathbb{Z}_+,
\end{equation}
is necessary and sufficient for the solvability of the moment problem~(\ref{f1_1}).

For a recent discussion on the truncated matrix
Hamburger moment problems we refer to the paper~\cite{cit_820_DFKMT} and
references therein.
It is worth mentioning that for the truncated moment problems much is done for the degenerate case, as well.
The case of the full moment problem~(\ref{f1_1}) is not such investigated. In~\cite{cit_850_Z}
we presented an analytic description of all solutions of the matrix Hamburger moment problem~(\ref{f1_1})
under condition~(\ref{f1_4}).
The main aim of our present investigation is to obtain a Nevanlinna-type formula for the moment
problem~(\ref{f1_1}) in a general case.
We only assume that condition~(\ref{f1_4}) holds and the moment problem~(\ref{f1_1}) is indeterminate
(but not necessarily completely indeterminate). We express the matrix coefficients of the corresponding
linear fractional transformation in terms of the given moments.
Some necessary and sufficient conditions for the determinacy of the moment problem~(\ref{f1_1})
in terms of the prescribed moments
are given.

\noindent
{\bf Notations.}  As usual, we denote by $\mathbb{R}, \mathbb{C}, \mathbb{N}, \mathbb{Z}, \mathbb{Z}_+$
the sets of real numbers, complex numbers, positive integers, integers, non-negative integers,
respectively; $\mathbb{C}_+ = \{ z\in \mathbb{C}:\ \mathop{\rm Im}\nolimits z > 0\}$,
$\mathbb{D} = \{ z\in \mathbb{C}:\ |z|<1 \}$, $\mathbb{T} = \{ z\in \mathbb{C}:\ |z|=1 \}$.
The notation $k\in\overline{0,\rho}$ means that $k\in \mathbb{Z}_+$, $k\leq\rho$, if $\rho <\infty$;
or $k\in \mathbb{Z}_+$, if $\rho = \infty$.
The set of all complex matrices of size $(m\times n)$ we denote by $\mathbb{C}_{m\times n}$, $m,n\in \mathbb{N}$.
If $M\in \mathbb{C}_{m\times n}$ then $M^T$ denotes the transpose of $M$, and
$M^*$ denotes the complex conjugate of $M$. The identity matrix from $\mathbb{C}_{n\times n}$
we denote by $I_n$, $n\in \mathbb{N}$; $I_\infty = (\delta_{k,l})_{k,l=0}^\infty$, $\delta_{k,l}$ is
Kronecker's delta.
If a set $S$ has a finite number of elements, then its number of elements we denote by
$\mathop{\rm card}\nolimits(S)$.
If a set $S$ has an infinite number of elements, then $\mathop{\rm card}\nolimits(S) := \infty$.

For a separable Hilbert space $H$ we denote by $(\cdot,\cdot)_H$ and $\| \cdot \|_H$ the scalar
product and the norm in $H$, respectively. The indices may be omitted in obvious cases.

\noindent
For a linear operator $A$ in $H$ we denote by $D(A)$ its domain, by $R(A)$ its range, and by
$A^*$ we denote its adjoint if it exists. If $A$ is invertible, then $A^{-1}$ means its inverse.
If $A$ is closable, then $\overline{A}$ means its closure.
If $A$ is bounded, then $\| A \|$ stands for its operator norm.
For a set of elements $\{ x_n \}_{n\in K}$ in $H$, we
denote by $\mathop{\rm Lin}\nolimits\{ x_n \}_{n\in K}$ and
$\mathop{\rm \overline{span}}\nolimits\{ x_n \}_{n\in K}$ the linear span and the closed
linear span in the norm of $H$, respectively. Here $K$ is an arbitrary set of indices.
For a set $M\subseteq H$, we denote by $\overline{M}$ the closure of $M$ in the norm of $H$.
By $E_H$ we denote the identity operator in $H$, i.e. $E_H x = x$, $x\in H$.
Let $H_1$ be a subspace of $H$. By $P_{H_1} = P_{H_1}^{H}$ we denote the operator
of the orthogonal projection on $H_1$ in $H$.

\noindent
If $A$ is symmetric, we set $R_z = R_z(A) = (A-zE_H)^{-1}$, $z\in \mathbb{C}\backslash \mathbb{R}$.
If $V$ is isometric, we set $\mathcal{R}_\zeta = \mathcal{R}_\zeta(V) = (E_H-\zeta V)^{-1}$,
$\zeta\in \mathbb{C}\backslash \mathbb{T}$.

\section{The matrix Hamburger moment problem: the determinacy and a Nevanlinna-type formula.}
Let the matrix Hamburger moment problem~(\ref{f1_1}) be given and condition~(\ref{f1_4}) hold.
Set
\begin{equation}
\label{f2_6}
\Gamma = (S_{k+l})_{k,l=0}^\infty = \left(
\begin{array}{ccccc} S_0 & S_1 & \ldots & S_n & \ldots\\
S_1 & S_2 & \ldots & S_{n+1} & \ldots\\
\vdots & \vdots & \ddots & \vdots & \ldots\\
S_n & S_{n+1} & \ldots & S_{2n} & \ldots\\
\vdots & \vdots & \vdots & \vdots & \ddots\end{array}
\right).
\end{equation}
The matrix $\Gamma$ is a semi-infinite block matrix. It may be viewed as a usual semi-infinite matrix, as well. Let
\begin{equation}
\label{f2_6_1}
\Gamma = (\Gamma_{n,m})_{n,m=0}^\infty,\qquad \Gamma_{n,m}\in \mathbb{C},
\end{equation}
and
$$ S_n = (s_n^{k,l})_{k,l=0}^{N-1},\qquad s_n^{k,l}\in \mathbb{C},\ n\in \mathbb{Z}_+. $$
Notice that
\begin{equation}
\label{f2_7}
\Gamma_{rN+j,tN+n} = s_{r+t}^{j,n},\qquad 0\leq j,n \leq N-1;\quad r,t\in \mathbb{Z}_+.
\end{equation}
We need here some constructions from~\cite{cit_850_Z}.
By Theorem~1 in~\cite{cit_850_Z} (and this construction is well known),
there exist a Hilbert space $H$, and a sequence $\{ x_n \}_{n=0}^\infty$ in $H$,
such that $\mathop{\rm \overline{span}}\nolimits\{ x_n \}_{n=0}^\infty = H$, and
\begin{equation}
\label{f2_9}
(x_n,x_m)_H = \Gamma_{n,m},\qquad n,m\in \mathbb{Z}_+.
\end{equation}
We choose an arbitrary such a space $H$ and a sequence $\{ x_n \}_{n=0}^\infty$ in $H$,
and {\it fix them in the rest of the paper}.

Set $L := \mathop{\rm Lin}\nolimits\{ x_n \}_{n=0}^\infty$, and consider the following operator with the domain $L$:
\begin{equation}
\label{f2_11}
A x = \sum_{k=0}^\infty \alpha_k x_{k+N},\qquad x\in L,\ x = \sum_{k=0}^\infty \alpha_k x_{k},\ \alpha_k\in \mathbb{C}.
\end{equation}
This operator is correctly defined and symmetric.

Let $\widehat A$ be an arbitrary self-adjoint extension of $A$ in a Hilbert space $\widehat H\supseteq H$.
Let $R_z(\widehat A) = (\widehat A - z E_{\widehat H})^{-1}$ be the resolvent of
$\widehat A$ and $\{ \widehat E_\lambda\}_{\lambda\in \mathbb{R}}$
be the orthogonal left-continuous resolution of unity of $\widehat A$. Recall that the operator-valued function
$\mathbf R_z = P_H^{\widehat H} R_z(\widehat A)$ is called a { \it generalized resolvent} of $A$,
$z\in \mathbb{C}\backslash \mathbb{R}$.
The function
$\mathbf E_\lambda = P_H^{\widehat H} \widehat E_\lambda$, $\lambda\in \mathbb{R}$, is a {\it spectral
function} of a symmetric operator $A$.
There exists a one-to-one correspondence between generalized resolvents and spectral functions.
It is given by the following relation~(\cite{cit_2000_AG}):
\begin{equation}
\label{f2_14}
(\mathbf R_z f,g)_H = \int_\mathbb{R} \frac{1}{\lambda - z} d( \mathbf E_\lambda f,g)_H,\qquad
f,g\in H,\ z\in \mathbb{C}\backslash \mathbb{R}.
\end{equation}
By Theorem~2 in~\cite{cit_850_Z}, all solutions of the moment problem~(\ref{f1_1}) have the following form:
\begin{equation}
\label{f2_30}
M(\lambda) = (m_{k,j} (\lambda))_{k,j=0}^{N-1},\quad
m_{k,j} (\lambda) = ( \mathbf E_\lambda x_k, x_j)_H,
\end{equation}
where $\mathbf E_\lambda$ is a spectral function of the operator $A$.
Moreover, the correspondence between all spectral functions of $A$ and all solutions
of the moment problem is bijective.

\noindent
By~(\ref{f2_14}) and~(\ref{f2_30}) we conclude that the formula
\begin{equation}
\label{f2_30_1}
\int_\mathbb{R} \frac{1}{\lambda - z} dm_{k,j} (\lambda) =
(\mathbf R_z x_k, x_j)_H,\quad 0\leq k,j\leq N-1,\quad z\in \mathbb{C}\backslash \mathbb{R},
\end{equation}
establishes a one-to-one correspondence between all generalized resolvents of $A$ and all solutions
of the moment problem~(\ref{f1_1}).

Let $B$ be a closed symmetric operator in the Hilbert space $H$, with the domain $D(B)$,
$\overline{D(B)} = H$.  Set $\Delta_B(\lambda) = (B- \lambda E_H) D(B)$,
and $N_\lambda = N_\lambda(B) = H\ominus \Delta_B(\lambda)$, $\lambda\in \mathbb{C}\backslash \mathbb{R}$.
Consider an
arbitrary bounded linear operator $C$, which maps $N_i$ into $N_{-i}$.
For
\begin{equation}
\label{f2_41}
g = f + C\psi - \psi,\qquad f\in D(B),\ \psi\in N_i,
\end{equation}
we set
\begin{equation}
\label{f2_42}
B_C g = Bf + i C \psi + i \psi.
\end{equation}
The operator $B_C$ is said to be a {\it quasiself-adjoint extension of the operator $B$, defined by
the operator $C$}.
By Theorem~4 in~\cite{cit_850_Z}, the following relation:
\begin{equation}
\label{f2_45}
\int_\mathbb{R} \frac{1}{x- \lambda } d m_{k,j} (x) = ( (A_{F(\lambda)} -
\lambda E_H)^{-1} x_k, x_j)_H,\qquad \lambda\in \mathbb{C}_+,
\end{equation}
establishes a bijective correspondence  between all solutions of the moment problem~(\ref{f1_1}) and all
analytic in $\mathbb{C}_+$ operator-valued functions $F(\lambda)$, which values are contractions
which map $N_i(\overline{A})$ into $N_{-i}(\overline{A})$. Here $A_{F(\lambda)}$ is the
quasiself-adjoint extension of $\overline{A}$ defined by $F(\lambda)$.

Set
$$ y_k^- := (A-iE_H) x_k = x_{k+N} - i x_k, $$
$$ y_k^+ := (A+iE_H) x_k = x_{k+N} + i x_k,\qquad k\in \mathbb{Z}_+; $$
\begin{equation}
\label{f2_45_1}
L^- :=  \mathop{\rm Lin}\nolimits\{ y_k^- \}_{k=0}^\infty = (A-iE_H)D(A),\quad
L^+ :=  \mathop{\rm Lin}\nolimits\{ y_k^+ \}_{k=0}^\infty = (A+iE_H)D(A),
\end{equation}
$$ H^- := \overline{L^-} = (\overline{A}-iE_H)D(\overline{A}),\quad
H^+ := \overline{L^+} = (\overline{A}+iE_H)D(\overline{A}). $$

Let us apply the Gram-Schmidt orthogonalization procedure to the sequence $\{ y_k^- \}_{k=0}^\infty$,
removing the linear dependent elements if they appear. We shall get a sequence
$\mathfrak{A}^- = \{ u_k^- \}_{k=0}^{\tau^- -1}$, $0\leq\tau^-\leq +\infty$.
The case $\tau^- = 0$ means that $y_k^- = 0$, $\forall k\in \mathbb{Z}_+$, and $\mathfrak{A}^-$ is
an empty set.

\noindent
In a similar manner,
we apply the Gram-Schmidt orthogonalization procedure to the sequence $\{ y_k^+ \}_{k=0}^\infty$,
and obtain a sequence
$\mathfrak{A}^+ = \{ u_k^+ \}_{k=0}^{\tau^+ -1}$, $0\leq\tau^+\leq +\infty$.
The case $\tau^+ = 0$ means that $y_k^+ = 0$, $\forall k\in \mathbb{Z}_+$, and $\mathfrak{A}^- = \emptyset$.

If not empty, the set $\mathfrak{A}^\pm$ forms an orthonormal basis in $H^\pm$, respectively.
Notice that, by the construction, each element $u_k^\pm$, $k\in \overline{0,\tau^\pm -1}$, is a linear combination
of $y_j^\pm$, $0\leq j\leq k$, respectively. Let
\begin{equation}
\label{f2_45_2}
u_k^\pm = \sum_{j=0}^k \xi_{k;j}^\pm y_j^\pm,\qquad \xi_{k;j}^\pm\in \mathbb{C},\quad k\in \overline{0,\tau^\pm -1}.
\end{equation}
Observe that by~(\ref{f2_9}) we may write
$$ (x_n,u_k^\pm)_H = \sum_{j=0}^k \overline{\xi_{k;j}^\pm} (x_n, y_j^\pm)_H
= \sum_{j=0}^k \overline{\xi_{k;j}^\pm} (x_n, x_{j+N} \pm i x_j)_H $$
\begin{equation}
\label{f2_45_3}
= \sum_{j=0}^k \overline{\xi_{k;j}^\pm} (\Gamma_{n,j+N} \pm i \Gamma_{n,j}),\quad n\in \mathbb{Z}_+,\
k\in\overline{0,\tau^\pm -1}.
\end{equation}
By representation~(\ref{f2_45_1}), the condition $\tau^- = 0$ ($\tau^+ = 0$) is equivalent to the
condition $D(A) = \{ 0 \}$, and therefore to the condition $H=\{ 0 \}$.
By~(\ref{f2_7}),(\ref{f2_9}), the condition $H=\{ 0 \}$ is equivalent to the condition
$S_n = 0$, $\forall n\in \mathbb{Z}_+$.

We emphasize that the numbers $\xi_{k;j}$ in~(\ref{f2_45_2}) can be computed explicitly by using
relations~(\ref{f2_7}),(\ref{f2_9}). Moreover, the processes of orthogonalization which appear
in this paper are based on the use of relations~(\ref{f2_7}),(\ref{f2_9}). In fact, any norm or
any scalar product which appear during orthogonalization is expressed in terms of the prescribed
moments.

\begin{thm}
\label{t2_1}
Let the matrix Hamburger moment problem~(\ref{f1_1}) be given and condition~(\ref{f1_4}),
with $\Gamma_n$ from~(\ref{f1_3}), be satisfied. Let the operator $A$ in the Hilbert space $H$
be constructed as in~(\ref{f2_11}).
The following conditions are equivalent:

\begin{itemize}

\item[{\rm (A)}]
The moment problem~(\ref{f1_1}) is determinate;

\item[{\rm (B)}]
One of the defect numbers of $A$ is equal to zero (or the both of them are zero);

\item[{\rm (C)}] $S_r = 0$, $\forall r\in \mathbb{Z}_+$, or, $\exists S_l\not=0$, $l\in \mathbb{Z}_+$,
and one of the following conditions holds (or the both of them hold):
{ \begin{itemize}

  \item[(a)]
For each $n$, $0 \leq n\leq N-1$, the following equality holds:
\begin{equation}
\label{f2_45_4}
\Gamma_{n,n} = \sum_{k=0}^{\tau^- -1}
\left|
\sum_{j=0}^k \overline{\xi_{k;j}^-} (\Gamma_{n,j+N} - i \Gamma_{n,j})
\right|^2;
\end{equation}

  \item[(b)]
For each $n$, $0 \leq n\leq N-1$, the following equality holds:
\begin{equation}
\label{f2_45_5}
\Gamma_{n,n} = \sum_{k=0}^{\tau^+ -1}
\left|
\sum_{j=0}^k \overline{\xi_{k;j}^+} (\Gamma_{n,j+N} + i \Gamma_{n,j})
\right|^2.
\end{equation}
Here $\Gamma_{\cdot,\cdot}$ are from~(\ref{f2_6_1}), and
$\xi_{\cdot,\cdot}^\pm$ are from~(\ref{f2_45_2}).

 \end{itemize} }

\end{itemize}

If the above conditions are satisfied then the unique solution of the moment problem~(\ref{f1_1})
is given by the following relation:
\begin{equation}
\label{f2_45_6}
M(t) = (m_{k,j} (t))_{k,j=0}^{N-1},\quad
m_{k,j} (t) = (E_t x_k, x_j)_H,
\end{equation}
where $E_t$ is the left-continuous orthogonal resolution of unity of the self-adjoint
operator $A$.
\end{thm}
{\bf Proof. }
(A)$\Rightarrow$(B).
If the both defect numbers are greater then zero, then we can choose unit vectors
$u_1\in N_i(\overline{A})$ and $u_2\in N_{-i}(\overline{A})$.
We set
$$ F(\lambda) (c u_1 + u) = c u_2,\qquad c\in \mathbb{C},\ u\in \Delta_{\overline A}(i). $$
On the other hand, we set
$\widetilde F(\lambda) \equiv 0$. Functions $F(\lambda)$ and $\widetilde F(\lambda)$ produce
different solutions of the moment problem~(\ref{f1_1}) by relation~(\ref{f2_45}).

\noindent
(B)$\Rightarrow$(A). If one of the defect numbers is zero, then the only admissible function $F(\lambda)$
in relation~(\ref{f2_45}) is $F(\lambda)\equiv 0$.

\noindent
(B)$\Rightarrow$(C).
If $H= \{ 0 \}$ then condition (C) holds. Let $H\not= \{ 0 \}$.

Notice that by~(\ref{f2_9}) and~(\ref{f2_45_3}), condition (C),(a) may be written as
$$ \| x_n \|^2 = \sum_{k=0}^{\tau^- -1} \left|
(x_n,u^-_k)_H
\right|^2,\qquad n=0,1,\ldots, N-1; $$
while condition (C),(b) is equivalent to
$$ \| x_n \|^2 = \sum_{k=0}^{\tau^+ -1} \left|
(x_n,u^+_k)_H
\right|^2,\qquad n=0,1,\ldots, N-1. $$
Therefore condition (C),(a) is equivalent to relations:
\begin{equation}
\label{f2_45_7}
x_n \in H^-,\qquad n=0,1,\ldots, N-1;
\end{equation}
and
condition (C),(b) is equivalent to condition:
\begin{equation}
\label{f2_45_8}
x_n \in H^+,\qquad n=0,1,\ldots, N-1.
\end{equation}
By the formula~(37) in~\cite[p. 278]{cit_850_Z}, each element of $L$ belongs to the linear span of
elements $\{ x_n \}_{n=0}^\infty$, $\{ y_k^- \}_{k=0}^\infty$, as well as to the linear span of
elements $\{ x_n \}_{n=0}^\infty$, $\{ y_k^+ \}_{k=0}^\infty$. Consequently,
condition~(\ref{f2_45_7})  is equivalent to the condition
\begin{equation}
\label{f2_45_9}
H = H^-,
\end{equation}
and condition~(\ref{f2_45_8})  is equivalent to the condition
\begin{equation}
\label{f2_45_10}
H = H^+.
\end{equation}
Since one of the defect numbers is equal to zero then either~(\ref{f2_45_9}), or~(\ref{f2_45_10})
holds.

\noindent
(C)$\Rightarrow$(B).
If $H= \{ 0 \}$ then condition (B) holds. Let $H\not= \{ 0 \}$.
If condition~(C),(a) (condition~(C),(b)) holds, then by the above considerations before~(\ref{f2_45_9}) we
obtain $H=H^-$ (respectively $H=H^+$). Therefore one of the defect numbers of $A$ is equal to zero.

The last assertion of the theorem follows from formula~(\ref{f2_30}).
$\Box$

We shall continue our considerations started before the statement of Theorem~\ref{t2_1}. In what follows
we assume that {\it the moment problem~(\ref{f1_1}) is indeterminate}. Let the defect numbers of
$A$ are equal to $\delta = \delta(A) = \dim H\ominus H^-$, and
$\omega = \omega(A) = \dim H\ominus H^+$,
$\delta,\omega\geq 1$.

For simplicity of notations we set $\tau := \tau^-$, and
$$ u_k := u_k^-,\qquad k\in\overline{0,\tau -1}. $$
Let us apply the Gram-Schmidt orthogonalization procedure to the vectors
$$ \{ u_k \}_{k=0}^{\tau -1}, \{ x_n \}_{n=0}^{N-1}. $$
Notice that the elements $\{ u_k \}_{k=0}^{\tau -1}$ are already orthonormal.
Then we get an orthonormal set in $H$:
$$ \mathfrak{A}_u := \{ u_k \}_{k=0}^{\tau -1} \cup \{ u_l' \}_{l=0}^{\delta-1}. $$
Notice that $\mathfrak{A}' := \{ u_l' \}_{l=0}^{\delta-1}$ is an orthonormal basis
in $H\ominus H^-$.

Set
\begin{equation}
\label{f2_45_11}
V = V_{\overline{A}} = (\overline{A} + iE_H)(\overline{A} - iE_H)^{-1} = E_H + 2i (\overline{A} - iE_H)^{-1}.
\end{equation}
The operator $V$ is a closed isometric operator with the domain $H^-$ and the range $H^+$.
Set
$$ v_k := V u_k,\qquad k\in \overline{0,\tau -1}. $$
Observe that by~(\ref{f2_45_2}) we may write
$$ v_k = \sum_{j=0}^k \xi_{k;j}^- V y_j^- =
\sum_{j=0}^k \xi_{k;j}^- y_j^+,\qquad k\in \overline{0,\tau -1}. $$
Notice that
$$ \mathfrak{A}_v^- := \{ v_k \}_{k=0}^{\tau -1}, $$
is an orthonormal basis in $H^+$.

\noindent
Let us apply the Gram-Schmidt orthogonalization procedure to the vectors
$$ \{ v_k \}_{k=0}^{\tau -1}, \{ x_n \}_{n=0}^{N-1}. $$
The elements $\{ v_k \}_{k=0}^{\tau -1}$ are already orthonormal.
Then we get another orthonormal basis in $H$:
$$ \mathfrak{A}_v := \{ v_k \}_{k=0}^{\tau -1} \cup \{ v_l' \}_{l=0}^{\omega-1}. $$
Observe that $\mathfrak{A}_v' := \{ v_l' \}_{l=0}^{\omega-1}$ is an orthonormal basis
in $H\ominus H^+$.

Let $\mathbf{R}_\lambda$ be an arbitrary generalized resolvent of the operator $A$. Let us check that
$$ (\mathbf{R}_z x_k, x_j)_H $$
$$ =
\frac{1}{z^2 +1} (\mathbf{R}_z y_k^-, y_j^-)_H
- \frac{1}{z^2 +1} (x_{k+N},x_j)_H
- \frac{z}{z^2 +1} (x_k,x_j)_H, $$
\begin{equation}
\label{f2_45_12}
 z\in \mathbb{C}_+\backslash\{ i \} ,\quad 0\leq k,j\leq N-1.
\end{equation}
In fact, let $\widetilde A\supseteq A$ be a self-adjoint operator in a Hilbert space
$\widetilde H\supseteq H$, such that $P^{\widetilde{H}}_H R_z(\widetilde A) = \mathbf{R}_z$,
$z\in \mathbb{C}\backslash \mathbb{R}$.
Then
$$ (\mathbf{R}_z x_k, x_j)_H = (R_z(\widetilde A) (A-iE_H)^{-1} (A-iE_H) x_k, x_j)_{\widetilde H} $$
$$ = (R_z(\widetilde A) R_i(\widetilde A) y_k^-, x_j)_{\widetilde H}
= \frac{1}{z - i} ((R_z(\widetilde A) - R_i(\widetilde A)) y_k^-, x_j)_{\widetilde H} $$
\begin{equation}
\label{f2_45_13}
 \frac{1}{z - i} (R_z(\widetilde A) y_k^-, x_j)_{\widetilde H} -
\frac{1}{z - i} (x_k, x_j)_{\widetilde H};
\end{equation}
$$ (R_z(\widetilde A) y_k^-, x_j)_{\widetilde H} =
(R_z(\widetilde A) y_k^-, R_i(\widetilde A) y_j^-)_{\widetilde H} $$
$$ = (R_{-i}(\widetilde A) R_z(\widetilde A) y_k^-, y_j^-)_{\widetilde H}
= -\frac{1}{i+z} ((R_{-i}(\widetilde A) - R_z(\widetilde A)) y_k^-, y_j^-)_{\widetilde H} $$
\begin{equation}
\label{f2_45_14}
 = -\frac{1}{i+z} (y_k^-, x_j)_{\widetilde H} +
\frac{1}{i+z} (\mathbf{R}_z y_k^-, y_j^-)_{\widetilde H}.
\end{equation}
By substitution~(\ref{f2_45_14}) into~(\ref{f2_45_13}), we get~(\ref{f2_45_12}).

Let $\widehat U\supseteq V$ be an arbitrary unitary extension of
$V$ in a Hilbert space $\widehat H\supseteq H$. Recall~\cite{cit_2100_Ch}
that the following function:
\begin{equation}
\label{f2_45_15}
\mathbf{R}_\zeta(V) =
P^{\widehat H}_H (E_{\widehat H} - \zeta \widehat U)^{-1},\qquad \zeta\in \mathbb{C}\backslash \mathbb{T},
\end{equation}
is said to be a {\it generalized resolvent} of $V$.

Observe that the generalized resolvents of $V$ and $\overline {A}$ are connected by the
following relation~\cite[pp. 370-371]{cit_2200_Ch}:
\begin{equation}
\label{f2_45_16}
(1-\zeta)\mathbf{R}_\zeta(V) =
E_H + (z -i) \mathbf{R}_z(\overline{A}),\qquad z\in \mathbb{C}_+,\ \zeta = \frac{z-i}{z+i} \in \mathbb{D}.
\end{equation}
(The latter relation follows from the fact that the usual resolvents of $V$
and $\overline{A}$ are related by a similar relation, and then one applies the projection operator
$P^{\widehat H}_H$ to the both sides of that relation.)
Correspondence~(\ref{f2_45_16}) between all generalized resolvents of $V$ and
all generalized resolvents of $\overline {A}$ is bijective.
Then
\begin{equation}
\label{f2_45_17}
\mathbf{R}_z(\overline{A}) =
\frac{2i}{z^2+1}\mathbf{R}_{\frac{z-i}{z+i}} (V) -
\frac{1}{z-i} E_H,\qquad z\in \mathbb{C}_+\backslash\{ i \}.
\end{equation}
By~(\ref{f2_45_17}),(\ref{f2_45_12}) and~(\ref{f2_9}) we get
$$ (\mathbf{R}_z x_k, x_j)_H
= \frac{2i}{(z^2 +1)^2} (\mathbf{R}_{\frac{z-i}{z+i}}(V_{ \overline{A} }) y_k^-, y_j^-)_H -
\frac{1}{(z^2+1)(z-i)}
\varphi_{j,k}(z), $$
\begin{equation}
\label{f2_45_18}
z\in \mathbb{C}_+\backslash\{ i \},\quad 0\leq k,j\leq N-1,
\end{equation}
where
$$ \varphi_{j,k}(z) :=
\Gamma_{k+N,j+N}-i\Gamma_{k+N,j}-i\Gamma_{k,j+N}+\Gamma_{k,j} + (z-i)\Gamma_{k+N,j} + z(z-i)\Gamma_{k,j} $$
\begin{equation}
\label{f2_45_19}
= \Gamma_{k+N,j+N} - i\Gamma_{k,j+N} + (z-2i)\Gamma_{k+N,j} + (z^2-iz+1)\Gamma_{k,j},\ z\in \mathbb{C}_+.
\end{equation}
Observe that an arbitrary generalized resolvent $\mathbf{R}_\zeta$ of the closed isometric operator $V_{\overline A}$
has the following representation~\cite[Theorem 3]{cit_2100_Ch}:
\begin{equation}
\label{f2_45_20}
\mathbf R_{\zeta} =
\left[
E - \zeta ( V \oplus \Phi_\zeta )
\right]^{-1},\qquad
\zeta\in \mathbb{D}.
\end{equation}
Here $\Phi_\zeta$ is an analytic in $\mathbb{D}$ operator-valued function which values are
linear contractions from $H\ominus H^-$ into $H\ominus H^+$.
The correspondence between all such functions $\Phi_\zeta$ and all generalized resolvents of $V$
is bijective.

\noindent
By~(\ref{f2_45_18}),(\ref{f2_45_20}) we get
$$ (\mathbf{R}_z x_k, x_j)_H =
\frac{2i}{(z^2 +1)^2} \left(\left[
E - \frac{z-i}{z+i} ( V \oplus \Phi_{\frac{z-i}{z+i}} )
\right]^{-1}
y_k^-, y_j^- \right)_H
$$
\begin{equation}
\label{f2_45_21}
- \frac{1}{(z^2+1)(z-i)}
\varphi_{j,k}(z),\qquad z\in \mathbb{C}_+\backslash\{ i \},\quad 0\leq k,j\leq N-1,
\end{equation}
where $\Phi_\cdot$ is an analytic in $\mathbb{D}$ operator-valued function which values are
linear contractions from $H\ominus H^-$ into $H\ominus H^+$.

\noindent
By~(\ref{f2_30_1}) and~(\ref{f2_45_21}) we conclude that the formula
$$ \int_\mathbb{R} \frac{1}{\lambda - z} dm_{k,j} (\lambda) $$
$$ =
\frac{2i}{(z^2 +1)^2} \left(\left[
E - \frac{z-i}{z+i} ( V \oplus \Phi_{\frac{z-i}{z+i}} )
\right]^{-1}
y_k^-, y_j^- \right)_H - \frac{1}{(z^2+1)(z-i)}
\varphi_{j,k}(z), $$
\begin{equation}
\label{f2_45_22}
0\leq k,j\leq N-1,\quad z\in \mathbb{C}_+\backslash \{ i  \},
\end{equation}
establishes a one-to-one correspondence between all
analytic in $\mathbb{D}$ operator-valued functions $\Phi_\cdot$, which values are
linear contractions from $H\ominus H^-$ into $H\ominus H^+$, and all solutions
$M(\lambda)=(m_{k,j}(\lambda))_{k,j=0}^{N-1}$
of the moment problem~(\ref{f1_1}).

It turns out that formula~(\ref{f2_45_22}) is more convenient then formula~(\ref{f2_45}), in order
to obtain a Nevanlinna-type formula for the moment problem~(\ref{f1_1}).

Denote by $\mathcal{M}_{1,\zeta}(\Phi)$ the matrix of the operator
$E_H - \zeta ( V \oplus \Phi_\zeta )$ in the basis $\mathfrak{A}_u$, $\zeta\in \mathbb{D}$. Here
$\Phi_\zeta$ is an analytic in $\mathbb{D}$ operator-valued function, which values are
linear contractions from $H\ominus H^-$ into $H\ominus H^+$. Then
$$ \mathcal{M}_{1,\zeta}(\Phi) =
\left(
\begin{array}{cc} A_{0,\zeta} & B_{0,\zeta}(\Phi) \\
C_{0,\zeta} & D_{0,\zeta}(\Phi) \end{array}
\right),
$$
where
$$ A_{0,\zeta} =
\left( \left( \left[
E_H - \zeta ( V \oplus \Phi_\zeta ) \right] u_k, u_j
\right)_H
\right)_{j,k=0}^{\tau-1}
=
\left( \left(
u_k - \zeta V u_k, u_j
\right)_H
\right)_{j,k=0}^{\tau-1} $$
\begin{equation}
\label{f2_50}
=
I_\tau - \zeta \left( \left(
v_k, u_j
\right)_H
\right)_{j,k=0}^{\tau-1},
\end{equation}
$$ B_{0,\zeta}(\Phi) =
\left( \left( \left[
E_H - \zeta ( V \oplus \Phi_\zeta ) \right] u_k', u_j
\right)_H
\right)_{0\leq j\leq \tau-1,\ 0\leq k\leq \delta-1} $$
$$ =
\left( \left(
u_k' - \zeta \Phi_\zeta u_k', u_j
\right)_H
\right)_{0\leq j\leq \tau-1,\ 0\leq k\leq \delta-1} $$
$$ =
-\zeta \left( \left(
\Phi_\zeta u_k', u_j
\right)_H
\right)_{0\leq j\leq \tau-1,\ 0\leq k\leq \delta-1}, $$
$$ C_{0,\zeta} =
\left( \left( \left[
E_H - \zeta ( V \oplus \Phi_\zeta ) \right] u_k, u_j'
\right)_H
\right)_{0\leq j\leq \delta-1,\ 0\leq k\leq \tau-1} $$
$$ =
\left( \left(
u_k - \zeta V u_k, u_j'
\right)_H
\right)_{0\leq j\leq \delta-1,\ 0\leq k\leq \tau-1} $$
\begin{equation}
\label{f2_51}
=
- \zeta \left( \left(
v_k, u_j'
\right)_H
\right)_{0\leq j\leq \delta-1,\ 0\leq k\leq \tau-1},
\end{equation}
$$ D_{0,\zeta}(\Phi) =
\left( \left( \left[
E_H - \zeta ( V \oplus \Phi_\zeta ) \right] u_k', u_j'
\right)_H
\right)_{0\leq j\leq \delta-1,\ 0\leq k\leq \delta-1} $$
$$ =
\left( \left(
u_k' - \zeta \Phi_\zeta u_k', u_j'
\right)_H
\right)_{0\leq j\leq \delta-1,\ 0\leq k\leq \delta-1} $$
$$ = I_\delta - \zeta
\left( \left(
\Phi_\zeta u_k', u_j'
\right)_H
\right)_{0\leq j\leq \delta-1,\ 0\leq k\leq \delta-1},\ \zeta\in \mathbb{D}. $$
Notice that matrices $A_{0,\zeta},C_{0,\zeta}$, $\zeta\in \mathbb{D}$, can be calculated explicitly using
relations~(\ref{f2_9}) and~(\ref{f2_7}).

Denote by $F_\zeta$, $\zeta\in \mathbb{D}$, the matrix of the operator $\Phi_\zeta$,
acting from $H\ominus H^-$ into $H\ominus H^+$, with respect to the bases
$\mathfrak{A}'$ and $\mathfrak{A}_v'$:
$$ F_\zeta = (f_\zeta(j,k))_{0\leq j\leq \omega - 1,\ 0\leq k\leq \delta-1},\qquad $$
$$ f_\zeta(j,k) := (\Phi_\zeta u_k', v_j')_H. $$
Then
$$ \Phi_\zeta u_k' = \sum_{l=0}^{\omega-1} f_\zeta(l,k) v_l',\quad 0\leq k\leq \delta-1, $$
and
$$ B_{0,\zeta}(\Phi) =
-\zeta \left( \left(
\sum_{l=0}^{\omega-1} f_\zeta(l,k) v_l', u_j
\right)_H
\right)_{0\leq j\leq \tau-1,\ 0\leq k\leq \delta-1} $$
$$ = -\zeta \left( \sum_{l=0}^{\omega-1} \left(
 v_l', u_j \right)_H f_\zeta(l,k)
\right)_{0\leq j\leq \tau-1,\ 0\leq k\leq \delta-1},\quad \zeta\in \mathbb{D}. $$
Set
\begin{equation}
\label{f2_52}
W := \left( \left(  v_l', u_j \right)_H
\right)_{0\leq j\leq \tau-1,\ 0\leq l\leq \omega-1}.
\end{equation}
Then
$$ B_{0,\zeta}(\Phi) = -\zeta W F_\zeta,\qquad \zeta\in \mathbb{D}. $$
We may write
$$ D_{0,\zeta}(\Phi) = I_\delta - \zeta
\left( \left(
\sum_{l=0}^{\omega-1} f_\zeta(l,k) v_l', u_j'
\right)_H
\right)_{0\leq j\leq\delta-1,\ 0\leq k\leq\delta-1} $$
$$ = I_\delta - \zeta
\left( \sum_{l=0}^{\omega-1}
\left( v_l', u_j' \right)_H f_\zeta(l,k)
\right)_{0\leq j\leq \delta-1,\ 0\leq k\leq \delta-1},\quad \zeta\in \mathbb{D}. $$
Set
\begin{equation}
\label{f2_53}
T :=
\left( \left( v_l', u_j' \right)_H
\right)_{0\leq j\leq \delta-1,\ 0\leq l\leq \omega-1}.
\end{equation}
Then
$$ D_{0,\zeta}(\Phi) = I_\delta - \zeta T F_\zeta,\qquad \zeta\in \mathbb{D}. $$
Thus, we may write
$$ \mathcal{M}_{1,\zeta}(\Phi) =
\left(
\begin{array}{cc} A_{0,\zeta} & -\zeta W F_\zeta \\
C_{0,\zeta} & I_\delta - \zeta T F_\zeta \end{array}
\right),\quad \zeta\in \mathbb{D},
$$
where $A_{0,\zeta}$, $C_{0,\zeta}$ are given by~(\ref{f2_50}),(\ref{f2_51}),
and $W,T$ are given by~(\ref{f2_52}),(\ref{f2_53}).

Consider the block representation of the operator $E_H - \zeta ( V \oplus \Phi_\zeta )$ with respect
to the decomposition $H^- \oplus (H\ominus H^-)$:
\begin{equation}
\label{f2_53_1}
E_H - \zeta ( V \oplus \Phi_\zeta ) =
\left(
\begin{array}{cc} \mathcal{A}_{0,\zeta} & \mathcal{B}_{0,\zeta}(\Phi) \\
\mathcal{C}_{0,\zeta} & \mathcal{D}_{0,\zeta}(\Phi) \end{array}
\right),\qquad \zeta\in \mathbb{D}.
\end{equation}
Of course, the matrices of operators $\mathcal{A}_{0,\zeta}$, $\mathcal{B}_{0,\zeta}$, $\mathcal{C}_{0,\zeta}$,
$\mathcal{D}_{0,\zeta}$ are matrices $A_{0,\zeta}$, $B_{0,\zeta}$, $C_{0,\zeta}$, $D_{0,\zeta}$,
respectively.
Observe that the matrix $A_{0,\zeta}$ is invertible, since
$\mathcal{A}_{0,\zeta} = P_{H^-} (E_H - \zeta V) P_{H^-} = E_{H^-} - \zeta P_{H^-} V P_{H^-}$, is
invertible, $\zeta\in \mathbb{D}$.
Set
\begin{equation}
\label{f2_53_2}
V_0 := P_{H^-} V P_{H^-}.
\end{equation}
The matrix of $V_0$ in the basis $\mathfrak{A}^-$ we denote by $\mathfrak{V}$:
\begin{equation}
\label{f2_53_3}
\mathfrak{V} = \left( \left(  v_k, u_j \right)_H \right)_{j,k=0}^{\tau-1}.
\end{equation}
Observe that using
definitions of $v_k$,$u_j$, the elements of the matrix $\mathfrak{V}$ can be calculated explicitly by
the prescribed moments.

We may write for the resolvent function of $V_0$:
\begin{equation}
\label{f2_53_4}
\mathcal{R}_\zeta(V_0) =
\mathcal{A}_{0,\zeta}^{-1} = E_{H^-} + \sum_{k=1}^\infty V^k \zeta^k,\qquad \zeta\in \mathbb{D}.
\end{equation}
Then for the corresponding matrices we may write:
\begin{equation}
\label{f2_53_5}
A_{0,\zeta}^{-1} = I_{\infty} + \sum_{k=1}^\infty \mathfrak{V}_k \zeta^k,\qquad \zeta\in \mathbb{D},
\end{equation}
where
\begin{equation}
\label{f2_53_6}
\mathfrak{V}_k := \mathfrak{V}^k,\qquad k\in \mathbb{Z}^+.
\end{equation}
By the convergence in~(\ref{f2_53_5}) we mean the convergence of the corresponding entries of matrices.

Observe that the Frobenius formula for the inverse of the block matrix (\cite[p. 59]{cit_7000_G}) is still valid for
the block representations of bounded operators as in~(\ref{f2_53_1}), if the following operator
$$ \mathcal{H}_\zeta := \mathcal{D}_{0,\zeta} - \mathcal{C}_{0,\zeta} \mathcal{A}_{0,\zeta}^{-1}
\mathcal{B}_{0,\zeta}, $$
has a bounded inverse. This can be verified by the direct multiplication of the corresponding
block representations.
Notice that in  our case $\mathcal{H}_\zeta$ has a bounded inverse.
In fact, we may write
$$ \left(
\begin{array}{cc} E_{H^-} & 0 \\
-\mathcal{C}_{0,\zeta} \mathcal{A}_{0,\zeta}^{-1} & E_{H\ominus H^-} \end{array}
\right)
\left(
\begin{array}{cc} \mathcal{A}_{0,\zeta} & \mathcal{B}_{0,\zeta}(\Phi) \\
\mathcal{C}_{0,\zeta} & \mathcal{D}_{0,\zeta}(\Phi) \end{array}
\right)
=
\left(
\begin{array}{cc} \mathcal{A}_{0,\zeta} & \mathcal{B}_{0,\zeta}(\Phi) \\
0 & \mathcal{H}_{\zeta} \end{array}
\right). $$
Observe that
$$ \left(
\begin{array}{cc} E_{H^-} & 0 \\
-\mathcal{C}_{0,\zeta} \mathcal{A}_{0,\zeta}^{-1} & E_{H\ominus H^-} \end{array}
\right)^{-1}
=
\left(
\begin{array}{cc} E_{H^-} & 0 \\
\mathcal{C}_{0,\zeta} \mathcal{A}_{0,\zeta}^{-1} & E_{H\ominus H^-} \end{array}
\right). $$
Therefore the operator
$\mathcal{Q} := \left(
\begin{array}{cc} \mathcal{A}_{0,\zeta} & \mathcal{B}_{0,\zeta}(\Phi) \\
0 & \mathcal{H}_{\zeta} \end{array}
\right)$ is invertible.

\noindent
Suppose that there exists $y\in H\ominus H^-$, $y\not= 0$, such that $\mathcal{H}_\zeta y = 0$.
Set $u := - \mathcal{A}_{0,\zeta}^{-1} \mathcal{B}_{0,\zeta} y$. Then
$$ \left(
\begin{array}{cc} \mathcal{A}_{0,\zeta} & \mathcal{B}_{0,\zeta}(\Phi) \\
0 & \mathcal{H}_{\zeta} \end{array}
\right)
\left(
\begin{array}{cc} u \\
y \end{array}
\right) = 0. $$
This contradicts to the invertibility of $\mathcal{Q}$.
Since $\mathcal{H}_\zeta^{-1}$ acts in the finite-dimensional space $H\ominus H^-$, it is bounded.

Applying the Frobenius formula we get
\begin{equation}
\label{f2_53_7}
(E_H - \zeta ( V \oplus \Phi_\zeta ))^{-1} =
\left(
\begin{array}{cc} \mathcal{A}_{0,\zeta}^{-1} + \mathcal{A}_{0,\zeta}^{-1} \mathcal{B}_{0,\zeta}
\mathcal{H}_{\zeta}^{-1} \mathcal{C}_{0,\zeta} \mathcal{A}_{0,\zeta}^{-1}
& \ast \\
\ast & \ast \end{array}
\right),\qquad \zeta\in \mathbb{D},
\end{equation}
where by stars $\ast$ we denote the blocks which are not of interest for us.

\noindent
Denote by $\mathcal{M}_{2,\zeta}(\Phi)$ the matrix of the operator
$(E_H - \zeta ( V \oplus \Phi_\zeta ))^{-1}$ in the basis $\mathfrak{A}_u$, $\zeta\in \mathbb{D}$.
Then
$$ \mathcal{M}_{2,\zeta}(\Phi) =
\left(
\begin{array}{cc} A_{0,\zeta}^{-1} + A_{0,\zeta}^{-1} B_{0,\zeta}
(D_{0,\zeta} - C_{0,\zeta} A_{0,\zeta}^{-1}
B_{0,\zeta})^{-1}
C_{0,\zeta} A_{0,\zeta}^{-1}
& \ast \\
\ast & \ast \end{array}
\right) $$
\begin{equation}
\label{f2_53_8}
=
\left(
\begin{array}{cc} A_{0,\zeta}^{-1} -\zeta A_{0,\zeta}^{-1} W F_\zeta
( I_\delta - \zeta T F_\zeta +\zeta  C_{0,\zeta} A_{0,\zeta}^{-1}
W F_\zeta)^{-1}
C_{0,\zeta} A_{0,\zeta}^{-1}
& \ast \\
\ast & \ast \end{array}
\right),\ \zeta\in \mathbb{D}.
\end{equation}
Let $\{ u_j \}_{j=0}^{\rho-1}$ be a set of elements which were obtained by the  Gram-Schmidt orthogonalization of
$\{ y_k^- \}_{k=0}^{N-1}$. Observe that $\rho \geq 1$. In the opposite case we
have $y_k^- = 0$, $0\leq k\leq N-1$. By~(\ref{f2_45_22}) we obtain that the moment problem~(\ref{f1_1})
is determinate, what contradicts to our assumptions. Set
$$ H_\rho^- := \mathop{\rm Lin}\nolimits\{ y_k^- \}_{k=0}^{N-1}
= \mathop{\rm Lin}\nolimits\{ u_j \}_{j=0}^{\rho-1}. $$
Consider the following operator:
$$ \mathcal{J}_\zeta := P^H_{H_\rho^-} (E_H - \zeta ( V \oplus \Phi_\zeta ))^{-1} P^H_{H_\rho^-},\quad
\zeta\in \mathbb{D}, $$
as an operator in the (finite-dimensional) Hilbert space $H_\rho^-$. Its matrix in the basis
$\{ u_j \}_{j=0}^{\rho-1}$
we denote by
$J_\zeta$. It is given by
$$ J_\zeta =
A_{1,\zeta} -\zeta A_{2,\zeta} W F_\zeta
( I_\delta - \zeta T F_\zeta +\zeta  C_{0,\zeta} A_{0,\zeta}^{-1}
W F_\zeta)^{-1}
C_{0,\zeta} A_{3,\zeta},\quad \zeta\in \mathbb{D}. $$
Here $A_{1,\zeta}$ is a matrix standing in the first $\rho$ rows and the first $\rho$ columns
of the matrix $A_{0,\zeta}^{-1}$;
$A_{2,\zeta}$ is a matrix standing in the first $\rho$ rows of the matrix $A_{0,\zeta}^{-1}$;
$A_{3,\zeta}$ is a matrix standing in the first $\rho$ columns
of the matrix $A_{0,\zeta}^{-1}$.

Consider the following operator from $\mathbb{C}^N$ to $H_\rho^-$:
$$ \mathcal{K} \sum_{n=0}^{N-1} c_n \vec e_n = \sum_{n=0}^{N-1} c_n y_n^-,\qquad c_n\in \mathbb{C}, $$
where $\vec e_n = (\delta_{n,0},\delta_{n,1},\ldots,\delta_{n,N-1})\in \mathbb{C}^N$.
Let $K$ be the matrix of $\mathcal{K}$ with respect to the orthonormal bases
$\{ \vec e_n \}_{n=0}^{N-1}$ and $\{ u_j \}_{j=0}^{\rho-1}$:
\begin{equation}
\label{f2_55_1}
K = \left( \left( \mathcal{K} \vec e_k, u_j \right)_H \right)_{0\leq j\leq \rho-1,\ 0\leq k\leq N-1}
= \left( \left(y_k^-, u_j \right)_H \right)_{0\leq j\leq \rho-1,\ 0\leq k\leq N-1}.
\end{equation}
By~(\ref{f2_45_22}) we may write that
$$ \int_\mathbb{R} \frac{1}{\lambda - z} dm_{k,j} (\lambda) $$
$$ = \frac{2i}{(z^2 +1)^2} \left( P^H_{H^-_\rho}\left[
E - \zeta ( V \oplus \Phi_\zeta )
\right]^{-1} P^H_{H^-_\rho}
\mathcal{K} \vec e_k, \mathcal{K} \vec e_j \right)_H - \frac{1}{(z^2+1)(z-i)}
\varphi_{j,k}(z)  $$
$$ = \frac{2i}{(z^2 +1)^2} \left( \mathcal{K}^* \mathcal{J}_\zeta
\mathcal{K} \vec e_k, \vec e_j \right)_{\mathbb{C}^N} -  \frac{1}{(z^2+1)(z-i)}
\varphi_{j,k}(z), $$
\begin{equation}
\label{f2_55_2}
0\leq k,j\leq N-1,\quad z\in \mathbb{C}_+\backslash\{ i \},\quad \zeta = \frac{z-i}{z+i},
\end{equation}
establishes a one-to-one correspondence between all
analytic in $\mathbb{D}$ operator-valued functions $\Phi_\cdot$, which values are
linear contractions from $H\ominus H^-$ into $H\ominus H^+$, and all solutions
$M(\lambda)=(m_{k,j}(\lambda))_{k,j=0}^{N-1}$
of the moment problem~(\ref{f1_1}).

\noindent
Observe that
$\left( \mathcal{K}^* \mathcal{J}_\zeta
\mathcal{K} \vec e_k, \vec e_j \right)_{\mathbb{C}^N}$
is the element in the $j$-th row and $k$-th column of the matrix $\mathcal{M}_{3,\zeta}$ of the operator
$\mathcal{J}_{1,\zeta} := \mathcal{K}^* \mathcal{J}_\zeta \mathcal{K}$ in the basis $\{ e_n \}_{n=0}^{N-1}$.
We may write
$$ \mathcal{M}_{3,\zeta} = K^* J_\zeta K $$
$$ =
K^* A_{1,\zeta} K - \zeta K^* A_{2,\zeta} W F_\zeta
( I_\delta + \zeta  (C_{0,\zeta} A_{0,\zeta}^{-1} W - T) F_\zeta)^{-1}
C_{0,\zeta} A_{3,\zeta} K,\quad \zeta\in \mathbb{D}. $$
Set
$$ \Delta(z) := (\varphi_{j,k}(z))_{j,k=0}^{N-1},\qquad z\in \mathbb{C}_+. $$
Then the following relation
$$ \int_\mathbb{R} \frac{1}{\lambda - z} dM^T (\lambda) = \frac{2i}{(z^2 +1)^2}
K^* A_{1,\zeta} K
 - \frac{1}{(z^2+1)(z-i)}
\Delta(z)  $$
$$ - \frac{2i}{(z^2 +1)^2} \zeta K^* A_{2,\zeta} W F_\zeta
( I_\delta + \zeta  (C_{0,\zeta} A_{0,\zeta}^{-1} W - T) F_\zeta)^{-1}
C_{0,\zeta} A_{3,\zeta} K, $$
\begin{equation}
\label{f2_55_3}
z\in \mathbb{C}_+\backslash\{ i \},\quad \zeta = \frac{z-i}{z+i},
\end{equation}
establishes a one-to-one correspondence between all
analytic in $\mathbb{D}$, $\mathbb{C}_{\omega\times\delta}$-valued functions $F_\zeta$, which values are
such that $F_\zeta^* F_\zeta\leq I_\delta$, and all solutions
$M(\lambda)$
of the moment problem~(\ref{f1_1}).
Set
$$ \mathbf{A}(z) = 2i K^* A_{1,\zeta} K - (z+i) \Delta(z), $$
$$ \mathbf{B}(z) = - 2i \zeta K^* A_{2,\zeta} W, $$
$$ \mathbf{C}(z) = \zeta  (C_{0,\zeta} A_{0,\zeta}^{-1} W - T), $$
\begin{equation}
\label{f2_55_4}
\mathbf{D}(z) = C_{0,\zeta} A_{3,\zeta} K,\quad
z\in \mathbb{C}_+\backslash\{ i \},\quad \zeta = \frac{z-i}{z+i}.
\end{equation}
Then the right-hand side of~(\ref{f2_55_3}) becomes
$$ \frac{1}{(z^2+1)^2} \mathbf{A}(z) + \frac{1}{(z^2+1)^2} \mathbf{B}(z) F_\zeta
( I_\delta +  \mathbf{C}(z) F_\zeta)^{-1} \mathbf{D}(z). $$

\begin{thm}
\label{t2_2}
Let the matrix Hamburger moment problem~(\ref{f1_1}) be given and condition~(\ref{f1_4}),
with $\Gamma_n$ from~(\ref{f1_3}), be satisfied. Suppose that the moment problem is
indeterminate.
All solutions of the moment problem~(\ref{f1_1}) can be obtained from the following relation:
$$ \int_\mathbb{R} \frac{1}{\lambda - z} dM^T (\lambda) $$
\begin{equation}
\label{f2_59}
= \frac{1}{(z^2+1)^2} \mathbf{A}(z) + \frac{1}{(z^2+1)^2} \mathbf{B}(z) \mathbf{F}(z)
( I_\delta +  \mathbf{C}(z) \mathbf{F}(z))^{-1} \mathbf{D}(z),\quad z\in \mathbb{C}_+\backslash\{ i \},
\end{equation}
where $\mathbf{A}(z)$, $\mathbf{B}(z)$, $\mathbf{C}(z)$, $\mathbf{D}(z)$ are analytic in
$\mathbb{C}_+$, matrix-valued functions defined by~(\ref{f2_55_4}),
with values in $\mathbb{C}_{N\times N}$, $\mathbb{C}_{N\times \delta}$,
$\mathbb{C}_{\delta\times \omega}$, $\mathbb{C}_{\delta\times N}$, respectively.
Here $\mathbf{F}(z)$ is an analytic in $\mathbb{C}_+$,
$\mathbb{C}_{\omega\times\delta}$-valued function which values are
such that $\mathbf{F}(z)^* \mathbf{F}(z) \leq I_\delta$, $\forall z\in \mathbb{C}_+$.
Conversely, each analytic in $\mathbb{C}_+$, $\mathbb{C}_{\omega\times\delta}$-valued function
such that $\mathbf{F}(z)^* \mathbf{F}(z) \leq I_\delta$, $\forall z\in \mathbb{C}_+$,
generates by relation~(\ref{f2_59}) a solution of the moment problem~(\ref{f1_1}).
The correspondence between all
analytic in $\mathbb{C}_+$, $\mathbb{C}_{\omega\times\delta}$-valued functions such that
$\mathbf{F}(z)^* \mathbf{F}(z) \leq I_\delta$, $\forall z\in \mathbb{C}_+$,
and all solutions of the moment problem~(\ref{f1_1}) is bijective.
\end{thm}
{\bf Proof. }
The proof follows from the preceding considerations.
$\Box$

\begin{center}
{\large\bf The Nevanlinna-type formula for the matrix Hamburger moment problem in a general case.}
\end{center}
\begin{center}
{\bf S.M. Zagorodnyuk}
\end{center}

In this paper we obtain a Nevanlinna-type formula for the matrix Hamburger moment
problem in a general case.
We only assume that the problem is solvable and has more that one solution.
We express the matrix coefficients of the corresponding
linear fractional transformation in terms of the prescribed moments.
Necessary and sufficient conditions for the determinacy of the moment problem are given.

}
\end{document}